\documentclass[reqno,a4paper,12pt]{amsart} 

\usepackage{amsmath,amscd,amsfonts,amssymb}
\usepackage{mathrsfs,dsfont}

\usepackage{bbm}

\allowdisplaybreaks

\numberwithin{equation}{section}
\numberwithin{figure}{section}

\newcommand\C{\mathbb{C}}
\newcommand{\M}{\mathcal{M}}
\newcommand{\A}{\mathcal{A}}
\newcommand{\B}{\mathcal{B}}
\newcommand{\K}{\mathcal{K}}
\newcommand\eps{\varepsilon}

\newcommand{\la}{\langle}
\newcommand{\ra}{\rangle}

\renewcommand\le{\leqslant}
\renewcommand\ge{\geqslant}

\newcommand{\sign}{\operatorname{sign}}
\newcommand{\tr}{\operatorname{tr}}

\theoremstyle{plain}
\newtheorem{thm}{Theorem}[section]
\newtheorem{lem}[thm]{Lemma}

\newtheorem*{claim*}{Claim}

\theoremstyle{definition}

\newtheorem*{definition*}{Definition}
\newtheorem*{remarks*}{Remarks}
\newtheorem*{remark*}{Remark}

\newenvironment{enumerate-alph}
{\begin{enumerate}
		\addtolength{\itemsep}{5pt}
		}
	{\end{enumerate}}

\newenvironment{enumerate-num}
{\begin{enumerate}
		\addtolength{\itemsep}{5pt}
		}
	{\end{enumerate}}

\newenvironment{enumerate-text}
{\begin{enumerate}
		\addtolength{\itemsep}{5pt}
		}
	{\end{enumerate}}

\begin{document}
	
	\title
	[Rescaling unconditional Schauder frames]
	{Rescaling of unconditional Schauder frames in Hilbert spaces and completely bounded maps}
	
	\author{Anton Tselishchev}
	\address{St. Petersburg Department of Steklov Mathematical Institute, Fontanka 27, St. Petersburg 191023, Russia}
	\email{celis\_anton@pdmi.ras.ru}

	\thanks{The work was supported by the Foundation for the Advancement of Theoretical
		Physics and Mathematics “BASIS”}

	\begin{abstract}
		We prove that if every element $u$ in a Hilbert space $H$ admits a representation as unconditionally convergent series $$u=\sum_{k=1}^\infty \la u, y_k\ra x_k,$$ then there exist nonzero scalars $\{\alpha_k\}_{k=1}^\infty$ such that both sequences $\{\alpha_k x_k\}_{k=1}^\infty$ and $\{\overline{\alpha}_k^{-1}y_k\}_{k=1}^\infty$ are frames. Our result has the following equivalent reformulation: if $\Phi:\ell^\infty\to B(H)$ is a bounded linear map such that for every element of the unit vector basis $e_k$ in $\ell^\infty$ the operator $\Phi(e_k)$ has rank one, then $\Phi$ is completely bounded.
	\end{abstract}

	\maketitle


	\section{Introduction}
	
	\subsection{}
	Let $H$ be a (complex, separable) Hilbert space. The following definition is very well known: the sequence $\{x_k\}_{k=1}^\infty$ is called a frame in $H$ if there exist constants $C, c > 0$ such that for every $u\in H$
	\begin{equation}\label{frame_def}
		c\|u\|^2\le \sum_{k=1}^\infty |\la u, x_k \ra|^2 \le C \|u\|^2.
	\end{equation}
	If the sequence $\{x_k\}_{k=1}^\infty$ satisfies only the upper bound of \eqref{frame_def} then it is called $C$-Bessel.
	
	One of the key properties of frames is the existence of a dual frame: if the sequence $\{x_k\}_{k=1}^\infty$ satisfies the condition~\eqref{frame_def}, then there exists a (usually non-unique) frame $\{y_k\}_{k=1}^\infty$ such that every $u\in H$ can be represented as the following unconditionally convergent series in $H$:
	\begin{equation}\label{Schframe_def}
		u = \sum_{k=1}^\infty \la u, y_k \ra x_k.
	\end{equation}
	We refer the reader to \cite[Chapter 8]{Hei11} for the basic theory of frames in Hilbert spaces.
	
	\subsection{}
	If the sequence $\{(x_k, y_k)\}_{k=1}^\infty$ is such that every $u\in H$ admits a representation as unconditionally convergent series \eqref{Schframe_def}, then $\{(x_k, y_k)\}_{k=1}^\infty$ is called an \emph{unconditional Schauder frame}. This notion can be easily generalized to arbitrary Banach spaces and it often serves as a natural analogue of frames in a Banach space setting.
	
	Therefore, the natural question arises: suppose that $\{(x_k, y_k)\}_{k=1}^\infty$ is an unconditional Schauder frame in a Hilbert space $H$; is it true that both $\{x_k\}_{k=1}^\infty$ and $\{y_k\}_{k=1}^\infty$ are frames? The answer to this question is, of course, ``no'': indeed, if $\{\alpha_k\}_{k=1}^\infty$ is an arbitrary sequence of nonzero complex numbers, then $\{(\alpha_k x_k, \overline{\alpha}_k^{-1}y_k) \}_{k=1}^\infty$ is still an unconditional Schauder frame, and we can take the numbers $\{\alpha_k\}_{k=1}^\infty$ so large that $\{\alpha_k x_k\}_{k=1}^\infty$ is not a frame. It is probably worth noting that if, for example, both sequences $\{x_k\}_{k=1}^\infty$ and $\{y_k\}_{k=1}^\infty$ are bounded from below (i.e., $\|x_k\|, \|y_k\|\ge c > 0$ for all $k$), then they are indeed frames, see e.g. \cite[Proposition 4.4]{LT25} or \cite[Theorem 3.13]{HLLL14}.
	
	But what happens in the general situation if we allow the rescaling discussed above: can we always choose the numbers  $\{\alpha_k\}_{k=1}^\infty$ so that both $\{\alpha_k x_k\}_{k=1}^\infty$ and $\{\overline{\alpha}_k^{-1}y_k \}_{k=1}^\infty$ are frames? This problem was stated as a conjecture, in a slightly different form (on the language of frame multipliers), in \cite{SB13} and \cite{BFPS24}, and some partial results were obtained in these papers. It was also addressed in the paper \cite{HLLL14}, where its connection with the theory of completely bounded maps between $C^*$-algebras was discovered. Moreover, in \cite[Theorem 5.1]{HLLL14} the authors claim that there exists a counterexample to this conjecture, however, the proof contains a mistake.
	
	Our goal is to provide a positive answer to this open problem. 
	
	\begin{thm}\label{thm:main}
		Suppose that $\{(x_k, y_k)\}_{k=1}^\infty$ is an unconditional Schauder frame in a Hilbert space $H$. Then there exists a sequence of nonzero complex numbers $\{\alpha_k\}_{k=1}^\infty\subset\C$ such that both $\{\alpha_k x_k\}_{k=1}^\infty$ and $\{\overline{\alpha}_k^{-1}y_k \}_{k=1}^\infty$ are frames.
	\end{thm}
	
	This theorem is a consequence of the following slightly more general result.
	
	\begin{thm}\label{thm:main_1}
		Suppose that $\{x_k\}_{k=1}^\infty$ and $\{y_k\}_{k=1}^\infty$ are sequences in a Hilbert space $H$ such that for every $u\in H$ the series
		\begin{equation}\label{eq:1.3}
			Tu = \sum_{k=1}^\infty \la u,y_k \ra x_k
		\end{equation}
		converges unconditionally. Then there exists a sequence of nonzero complex numbers $\{\alpha_k\}_{k=1}^\infty$ such that both $\{\alpha_k x_k\}_{k=1}^\infty$ and $\{\overline{\alpha}_k^{-1}y_k \}_{k=1}^\infty$ are Bessel.
	\end{thm}
	
	It is also reasonable to ask what Bessel bound we can get for the sequences $\{\alpha_k x_k\}_{k=1}^\infty$ and $\{\overline{\alpha}_k^{-1}y_k\}$ --- we address this question in Theorem~\ref{thm:main_2} below.
	
	Let us briefly explain why Theorem~\ref{thm:main_1} implies Theorem~\ref{thm:main}. Since $\{(x_k, y_k)\}_{k=1}^\infty$ is an unconditional Schauder frame, the operator $T$ given by the formula \eqref{eq:1.3} is simply an identity operator. Theorem~\ref{thm:main_1} then gives us the complex numbers $\{\alpha_k\}_{k=1}^\infty$ such that the sequences $\{\alpha_k x_k\}_{k=1}^\infty$ and $\{\overline{\alpha}_k^{-1}y_k \}_{k=1}^\infty$ are Bessel. But since $\{(\alpha_k x_k, \overline{\alpha}_k^{-1}y_k) \}_{k=1}^\infty$ is also an unconditional Schauder frame, then these sequences must in fact be frames, see e.g. \cite[Theorem 3.13]{HLLL14} and also the proof of \cite[Proposition 4.4]{LT25} for a short self-contained presentation of this fact.
	
	Our proof of Theorem~\ref{thm:main_1} requires some non-trivial results from operator theory (about completely bounded maps) and a certain computation involving Rademacher functions which resembles the proof of the so-called ``little'' Grothendieck theorem.
	
	\subsection*{Acknowledgments} The author is grateful to Daniel Freeman for communicating this problem and to Nir Lev and Sergey Kislyakov for helpful remarks and discussions.
	
	\section{Preliminaries and reformulations of the main theorem}
	
	\subsection{} As we mentioned above, we will use the Rademacher functions in the proof. They can be defined by the formula $r_n(t) = \sign \sin (2^n \pi t)$ for all $n\ge 1$. We will use the Khintchine's inequality for $L^1$: there exists an absolute constant $\gamma >0$ such that for every finite sequence of complex numbers $\{a_k\}_{k=1}^n$
	\begin{equation}
		\int_0^1  \Big| \sum_{k=1}^n a_k r_k(t) \Big|\, dt\ge \gamma \cdot \Big( \sum_{k=1}^n |a_k|^2 \Big)^{1/2}.
	\end{equation} 
	In fact, it is known that the optimal value of $\gamma$ in this inequality is $1/\sqrt{2}$, see \cite{Sza76}.
	
	\subsection{}
	
	We will also need the notion of a trace norm on the space of $m\times m$ matrices. Recall that the trace-norm of an $m\times m$ matrix $B$ is defined as $\|B\|_1 = \tr |B|$, where $|B|=\sqrt{B^*B}$. The following simple lemma shows how to compute the trace-norm of a rank one matrix.
	
	\begin{lem}\label{lem:trace_norm_comput}
		Suppose that $B=(\alpha_i\beta_j)_{i,j=1}^m$ is a rank-one matrix. Then its trace-norm can be computed as follows:
		\begin{equation}\label{eq:trace_norm_comput}
			\|B\|_1=\tr|B| = \Big(\sum_{i=1}^m |\alpha_i|^2\Big)^{1/2}\cdot \Big(\sum_{j=1}^m |\beta_j|^2\Big)^{1/2}.
		\end{equation}
	\end{lem}
	
	The proof of this lemma is essentially a simple computation, we outline it below for completeness.
	
	\begin{proof}
		Put $\alpha=(\alpha_i)_{i=1}^m, \beta = (\overline{\beta}_j)_{j=1}^m\in\C^m$. Then the matrix $B$ is naturally identified with an operator $\alpha\otimes\beta$ on $\C^m$ which is defined as $(\alpha\otimes\beta)(x) = \la x, \beta\ra \alpha$. The matrix $B^*$ then corresponds to the operator $\beta\otimes\alpha$ and for every $x\in\C^m$ we have $(B^*B)(x) = B^*(\la x,\beta\ra\alpha)=\|\alpha\|^2 \la x, \beta \ra \beta$. It can be easily checked that then $|B|$ is the matrix associated with the operator defined as $|B| x = \frac{\|\alpha\|}{\|\beta\|}\la x,\beta\ra \beta$, and this operator has a unique nonzero eigenvalue $\|\alpha\|\cdot \|\beta\|$ (the corresponding eigenvector is $\beta$), so the identity \eqref{eq:trace_norm_comput} follows.
	\end{proof}
	
	For more information about the trace and the trace norm (in a more general context of bounded operators on a Hilbert space) see e.g. \cite[Section 3.4]{Ped89}.
	
	\subsection{}
	It was shown in \cite{BFPS24} that Theorem~\ref{thm:main_1} has the following equivalent reformulation.
	
	\begin{thm}\label{thm:main_2}
		There exists an absolute constant $K > 0$ satisfying the following property. Suppose that $\{x_k\}_{k=1}^n$ and $\{y_k\}_{k=1}^n$ are finite sequences of vectors in a finite-dimensional Hilbert space $H$ such that
		\begin{equation}\label{eq:2.2}
			\Big\| \sum_{k=1}^n \eps_k \la u, y_k \ra x_k \Big\|\le C\|u\|\quad \text{for all}\ u\in H \text{and}\ |\eps_k|\le 1.
		\end{equation}
		Then there exists a sequence of nonzero complex numbers $\{\alpha_k\}_{k=1}^n$ such that both $\{\alpha_k x_k\}_{k=1}^n$ and $\{\overline{\alpha}_k^{-1}y_k \}_{k=1}^n$ are $CK$-Bessel.
	\end{thm}
	
	Our goal is to prove this equivalent reformulation of Theorem~\ref{thm:main_1}. Obviously, we can assume without loss of generality that all vectors $\{x_k\}_{k=1}^\infty$ and $\{y_k\}_{k=1}^\infty$ are nonzero, and in what follows we will always work under this assumption.
	
	\subsection{} 
	It was noted in \cite{HLLL14} that Theorems~\ref{thm:main_1} and~\ref{thm:main_2} can be reformulated using the notion of completely bounded maps between $C^*$-algebras. Here we present the necessary definitions which are related to our problem and refer the reader to the book \cite{Pau02} for the general theory of completely bounded maps.
	
	We denote by $B(H)$ the algebra of all bounded operators on a Hilbert space $H$. Suppose that $\A\subset B(H)$ is a $C^*$-algebra. Then we can define the vector space $\M_m(\A)$ of $m\times m$-matrices whose entries are the elements of $\A$. It can be naturally identified with a space of operators on $H^m = \underbrace{H\oplus H\oplus \ldots \oplus H}_{m\, \text{times}}$. If we consider an operator norm on $\M_m(\A)$, then it becomes a Banach space (in fact, it also becomes a $C^*$-algebra).
	
	Now consider two $C^*$-algebras $\A$ and $\B$ (they can be viewed as subalgebras of $B(H_1)$ and $B(H_2)$ for two different Hilbert spaces $H_1$ and $H_2$). Suppose that $\Psi:A\to B$ is a bounded linear map. Then for every $m\ge 1$ we can define the linear map $\Psi_{(m)}:\M_m(\A)\to \M_m(\B)$ which simply acts entry-wise, i.e.,
	\begin{equation}
		\Psi_{(m)}((a_{i,j})_{i,j=1}^m) = (\Psi(a_{i,j}))_{i,j=1}^m.
	\end{equation}
	
	The map $\Psi$ is called \emph{completely bounded} if all maps $\Psi_{(m)}$ are uniformly bounded. If $\Psi$ is a completely bounded map, then we put $\|\Psi\|_{cb}=\sup_m \|\Psi_{(m)}\|$.
	
	\subsection{}
	If we are given the sequences $\{x_k\}_{k=1}^n$ and $\{y_k\}_{k=1}^n$ satisfying the condition \eqref{eq:2.2}, then we can consider the linear map $\Phi: \ell^\infty_n\to B(H)$ defined by
	\begin{equation}\label{eq:2.4}
		\Phi(a)u = \sum_{k=1}^n a(k)\la u, y_k\ra x_k
	\end{equation}
	for every $a = (a(k))_{k=1}^n\in\ell^\infty_n$. The inequality~\eqref{eq:2.2} then can be rewritten as 
	\begin{equation}
		\|\Phi(a)u\|_H\le C\|u\|_H\|a\|_{\ell^\infty_n},
	\end{equation}
	which in turn simply means that $\|\Phi\|\le C$. We also note that the map $\Phi$ does not depend on our rescaling, i.e., the sequences $\{\alpha_k x_k\}_{k=1}^n$ and $\{\overline{\alpha}_k^{-1}y_k \}_{k=1}^n$ give rise to the same map $\Phi$.
	
	Recall that $\ell^\infty_n$ is a $C^*$-algebra which can be identified with an algebra of diagonal operators on $\C^n$.
	
	It turns out that the conclusion of Theorem~\ref{thm:main_2} is closely related to the estimate of the norm $\|\Phi\|_{cb}$. Let us denote by $\{e_k\}_{k=1}^n$ the unit vector basis in $\ell^\infty_n$. Clearly, a map $\Phi:\ell^\infty_n\to B(H)$ has the form \eqref{eq:2.4} if and only if $\Phi(e_k)$ is a rank one operator for every $1\le k\le n$, where $\{e_k\}_{k=1}^n$ is a unit vector basis in $\ell^\infty_n$.
	
	\begin{thm}\label{thm:main_cb}
		There exists a positive constant $K > 0$ satisfying the following property. If $H$ is a finite-dimensional Hilbert space and $\Phi: \ell^\infty_n\to B(H)$ is a bounded linear map such that  $\Phi(e_k)$ is a rank one operator for every $1\le k\le n$, then $\|\Phi\|_{cb}\le K\|\Phi\|$.
	\end{thm}
	
	The fact that the problem about rescaling of unconditional Schauder frames into a dual pair of frames has an equivalent reformulation on the language of completely bounded maps was established in \cite[Section 3]{HLLL14}. We include a short proof that Theorem~\ref{thm:main_cb} implies Theorem~\ref{thm:main_2} in the final section of the present paper, which also shows that the validity of Theorem~\ref{thm:main_cb} for some constant $K$ implies that Theorem~$\ref{thm:main_2}$ (strictly speaking, for nonzero vectors $\{x_k\}_{k=1}^\infty$ and $\{y_k\}_{k=1}^\infty$) holds with the same constant. 
	
	\subsection{Remarks}
	
	1. The evaluation of the best possible constant $K$ in Theorems~\ref{thm:main_2} and~\ref{thm:main_cb} is beyond the scope of the present paper, although it seems to be an interesting problem. We only mention that by using the analogue of Khintchine's inequality for independent Gaussian random variables (as it is done in \cite[Section 5]{Pis12}) one can get the better constant $K$ than $K=2$ which is obtained in our proof below.
	
	2. Theorem~\ref{thm:main_cb} can be generalized in a following way: for every $r\ge 1$ there exists a constant $K(r)$ such that if $\mathrm{rank}\,\Phi(e_k)\le r$ for every $1\le k\le n$, then $\|\Phi\|_{cb}\le K(r) \|\Phi\|$. The proof of this fact can be obtained by an appropriate modification of the proof of Theorem~\ref{thm:main_cb} below. It would be interesting to establish the asymptotic behavior of the optimal constants $K(r)$, however we do not address this question in the present paper.
	
	3. It is known that Theorem~\ref{thm:main_cb} is not true for arbitrary linear maps $\Phi:\ell^\infty_n\to B(H)$. Indeed, for every $n\ge 1$ an example of a map $\Phi: \ell^\infty_n\to B(\C^{2^n})$ such that $\|\Phi\|_{cb}\ge \|\Phi_{(2^n)}\|\ge\sqrt{n} \|\Phi\|$ was constructed in \cite{Loe76}.
	
	4. Theorem~\ref{thm:main_cb} has the following infinite-dimensional version: if $\Phi:\ell^\infty\to B(H)$ is a bounded linear map such that $\Phi(e_k)$ is a rank one operator for every $k$, then $\Phi$ is completely bounded and $\|\Phi\|_{cb}\le K \|\Phi\|$ for some absolute constant $K$. In this formulation $H$, of course, denotes an inifinte-dimensional separable Hilbert space.

	\medskip
	
	In the next section we prove Theorem~\ref{thm:main_cb}.
	
	\section{The proof of Theorem~\ref{thm:main_cb}}
	
	\subsection{}
	We start by rewriting the inequality $\|\Phi\|_{cb}\le K\|\Phi\|$ which we have to prove.
	
	Since $\Phi(e_k)$ is a rank one operator, $\Phi$ is given by the formula \eqref{eq:2.4} for some $\{x_k\}_{k=1}^n$ and $\{y_k\}_{k=1}^n$. We need to show that $\|\Phi_{(m)}\|\le K \|\Phi\|$ for every $m\ge 1$, where $K$ is some absolute constant.
	
	Consider an element $A\in \M_m(\ell^\infty_n)$, $A=(a_{i,j}(k))_{i,j=1}^m$. Then $A$ can be viewed as a sequence of $m\times m$ matrices $\{A_k\}_{k=1}^n$, where $A_k=(a_{i,j}(k))_{i,j=1}^m$ are matrices whose entries are complex numbers. Then the norm of $A$ in the space $\M_m(\ell^\infty_n)$ can be computed as
	\begin{equation}\label{eq:3.2}
		\|A\|=\max_{1\le k\le n}\|A_k\|,
	\end{equation}
	where the norm of the matrix $A_k$ in the right-hand side stands for its norm as an operator on the space $\ell^2_n$. 
	
	We need to show that 
	\begin{equation}
		\|\Phi_{(m)}A\|\le K\cdot \|\Phi\|\cdot\|A\|.
	\end{equation} Recall that $\Phi_{(m)}(A)=(\Phi(a_{i,j}))_{i,j=1}^m\in\M_m(B(H))$ is viewed as an operator on the Hilbert space $H^m$. Therefore, in order to prove~\eqref{eq:3.2}, we need to verify the following inequality for arbitrary vectors $u_1, u_2, \ldots u_m\in H$:
	\begin{equation}
		\Big(\sum_{i=1}^m \Big\| \sum_{j=1}^m \Phi(a_{ij})u_j \Big\|^2\Big)^{1/2}\le K \cdot \|\Phi\|\cdot \|A\|\cdot \Big( \sum_{j=1}^m \|u_j\|^2 \Big)^{1/2}.
	\end{equation}
	
	Expanding $\Phi(a_{ij})u_j$ by the formula \eqref{eq:2.4}, we rewrite our inequality as
	
	\begin{equation}
		\Big(\sum_{i=1}^m \Big\| \sum_{j=1}^m \sum_{k=1}^n a_{i,j}(k)\la u_j, y_k\ra x_k \Big\|^2\Big)^{1/2}\le K \cdot \|\Phi\|\cdot \|A\|\cdot \Big( \sum_{j=1}^m \|u_j\|^2 \Big)^{1/2}.
	\end{equation}
	
	Using duality and changing the order of summation, we see that this inequality is equivalent to the estimate
	\begin{equation}\label{eq:key_to_prove}
		\Big|\sum_{k=1}^n \sum_{i=1}^m \sum_{j=1}^m a_{i,j}(k) \la u_j, y_k\ra \la x_k, v_i \ra\Big| \le  K \cdot \|\Phi\|\cdot \|A\|\cdot \Big( \sum_{j=1}^m \|u_j\|^2 \Big)^{1/2}\Big( \sum_{i=1}^m \|v_i\|^2 \Big)^{1/2}.
	\end{equation}
	Summing up, it only remains to prove the inequality \eqref{eq:key_to_prove} for arbitrary $\{u_j\}_{j=1}^m, \{v_i\}_{i=1}^m\subset H$.
	
	\subsection{}
	
	Recall that $\Phi$ is defined by the formula 
	\begin{equation}
		\Phi(a)u = \sum_{k=1}^n a(k)\la u, y_k\ra x_k
	\end{equation}
	and, by the definition of the quantity $\|\Phi\|$,
	\begin{equation}
		\|\Phi(a) u\|\le \|\Phi\|\cdot \|a\|_{\ell^\infty_n}\cdot \|u\|.
	\end{equation}
	It follows from the above two formulas that for arbitrary $u, v\in H$ we have
	\begin{equation}
		|\la \Phi(a)u, v \ra| = \Big| \sum_{k=1}^n a(k) \la u, y_k\ra \la x_k, v\ra \Big|\le \|\Phi\|\cdot \|a\|_{\ell^\infty_{n}}\cdot\|u\|\cdot\|v\|.
	\end{equation}
	By choosing an appropriate unimodular sequence $\{a(k)\}_{k=1}^n$ (depending on $u$ and $v$) we get the following inequality for all $u,v\in H$:
	\begin{equation}\label{eq:key_simple}
		\sum_{k=1}^n |\la u, y_k \ra|\cdot |\la v, x_k\ra|\le \|\Phi\|\cdot\|u\| \cdot\|v\|.
	\end{equation}
	
	\subsection{}
	
	Our next step is to apply Khintchine's inequality. The computations below are somewhat similar to the computations in the proof of the ``little'' Grothendieck theorem presented in \cite[Section 5]{Pis12}.
	
	Consider the Rademacher functions $\{r_j\}_{j=1}^m$ and for arbitrary $t, s\in [0,1]$ put $u=\sum_{j=1}^m r_j(t)u_j$ and $v = \sum_{i=1}^m r_i(s) v_i$ into the inequality~\eqref{eq:key_simple}. We get
	\begin{equation}
		\sum_{k=1}^n \Big| \la \sum_{j=1}^m r_j(t) u_j, y_k \ra \Big|\cdot \Big| \la \sum_{i=1}^m r_i(s)v_i, x_k \ra \Big|\le \|\Phi\|\cdot \Big\| \sum_{j=1}^m r_j(t) u_j \Big\|\cdot \Big\| \sum_{i=1}^m r_i(s)v_i \Big\|.
	\end{equation}
	Now we integrate this inequality with respect to $t$ and $s$ and write the following chain of inequalities:
	\begin{align}
		&\sum_{k=1}^n \Big( \int_0^1 \Big| \sum_{j=1}^m r_j(t)\la u_j,y_k\ra \Big|\, dt \Big)\cdot \Big( \int_0^1 \Big| \sum_{i=1}^m r_i(s)\la v_j,x_k\ra \Big|\, ds \Big)\label{eq:3.11}\\
		&\le \|\Phi\|\cdot \Big(\int_0^1 \Big\| \sum_{j=1}^m r_j(t) u_j \Big\|\, dt\Big) \cdot \Big(\int_0^1\Big\| \sum_{i=1}^m r_i(s)v_i \Big\|\, ds\Big)\\
		&\le  \|\Phi\|\cdot \Big(\int_0^1 \Big\| \sum_{j=1}^m r_j(t) u_j \Big\|^2\, dt\Big)^{1/2} \cdot \Big(\int_0^1\Big\| \sum_{i=1}^m r_i(s)v_i \Big\|^2\, ds\Big)^{1/2}\\
		&=\|\Phi\|\cdot \Big( \sum_{j=1}^m \|u_j\|^2 \Big)^{1/2}\cdot \Big(\sum_{i=1}^m \|v_i\|^2 \Big)^{1/2}.
	\end{align}
	The last equality can be easily verified using the fact that the Rademacher functions are pairwise orthogonal in $L^2$.
	
	Now we estimate the quantity \eqref{eq:3.11} from below using Khintchine's inequality: it is not smaller than
	\begin{equation}
		\gamma^2\cdot \sum_{k=1}^n \Big\{ \Big( \sum_{j=1}^m |\la u_j, y_k \ra|^2 \Big)^{1/2}\cdot \Big( \sum_{i=1}^m |\la v_i, x_k\ra|^2 \Big)^{1/2}\Big\}.
	\end{equation} 
	Summing up, we arrive at the following key estimate:
	\begin{equation}\label{eq:super_key_est}
		\begin{gathered}
			\sum_{k=1}^n \Big\{ \Big( \sum_{j=1}^m |\la u_j, y_k \ra|^2 \Big)^{1/2}\cdot \Big( \sum_{i=1}^m |\la v_i, x_k\ra|^2 \Big)^{1/2}\Big\}\\
			\le \gamma^{-2}\cdot \|\Phi\|\cdot \Big( \sum_{j=1}^m \|u_j\|^2 \Big)^{1/2}\cdot \Big(\sum_{i=1}^m \|v_i\|^2 \Big)^{1/2}.
		\end{gathered}
	\end{equation}
	
	\subsection{}
	
	It remains to show that the inequality \eqref{eq:super_key_est} implies \eqref{eq:key_to_prove}.
	
	For every $k$ consider the $m\times m$ matrix $B_k=(\la u_j, y_k\ra\cdot \la x_k, v_i\ra)_{i,j=1}^m$. Clearly, it is a rank one matrix. Then Lemma~\ref{lem:trace_norm_comput} shows that the left-hand side of the estimate~\eqref{eq:super_key_est} is equal to $\sum_{k=1}^n\|B_k\|_1$ (recall that $\|B_k\|_1 = \tr |B_k|$ is the trace-norm of $B_k$).
	
	The left-hand side of the inequality~\eqref{eq:key_to_prove} can in turn be rewritten as
	\begin{equation}
		\Big| \sum_{k=1}^n \tr(A_kB_k^t) \Big|,
	\end{equation}	
	where $B_k^t$ stands for a transpose of the matrix $B_k$. Therefore, the left-hand side of~\eqref{eq:key_to_prove} does not exceed
	\begin{align}
		&\sum_{k=1}^n |\tr(A_k B_k^t)|\le \sum_{k=1}^n \|A_k\|\cdot\|B_k\|_1\\
		&\le (\max_{1\le k\le n}\|A_k\|)\cdot \Big(\sum_{k=1}^n \|B_k\|_1\Big) = \|A\|\cdot \Big(\sum_{k=1}^n \|B_k\|_1\Big).
	\end{align}
	Here we used the well-known inequality $|\tr(AB)|\le \|A\|\cdot\|B\|_1$, see e.g. \cite[Lemma 3.4.10]{Ped89}.
	
	It remains to use the estimate of $\sum_{k=1}^n \|B_k\|_1$ given by~\eqref{eq:super_key_est}, and we are done: Theorem~\ref{thm:main_cb} is proved with constant $K=\gamma^{-2}=2$.
	
	\section{Relation between completely bounded maps and rescaling}
	
	\subsection{}
	
	Now we prove that Theorem~\ref{thm:main_cb} implies Theorem~\ref{thm:main_2}. Our proof relies on the generalization of Stinespring's representation for completely bounded maps, see~\cite[Theorem 8.4]{Pau02}. We formulate this result here.
	
	\begin{lem}\label{lem:Stinespring}
		Let $\A$ be a $C^*$-algebra with a unit and let $\Psi:\A\to B(H)$ be a completely bounded map with $\|\Psi\|_{cb}=1$. Then there exists a Hilbert space $\mathcal{K}$, a $*$-homomorphism $\pi:\A\to B(\mathcal{K})$ and isometries $V_{1,2}:H\to\mathcal{K}$ such that
		\begin{equation}
			\Psi(a)=V_1^* \pi(a)V_2
		\end{equation}
		for all $a\in\A$.
	\end{lem}
	
	\subsection{}
	
	Suppose that we are given the sequences $\{x_k\}_{k=1}^n$ and $\{y_k\}_{k=1}^n$ in a finite-dimensional Hilbert space $H$ which satisfy the conditions of Theorem~\ref{thm:main_2}. Then we can consider the linear map $\Phi:\ell^\infty_n\to B(H)$ defined by the formula~\eqref{eq:2.4}, i.e.
	\begin{equation}\label{eq:4.2}
		\Phi(a)u = \sum_{k=1}^n a(k)\la u, y_k\ra x_k.
	\end{equation}
	As we have already mentioned in Section~2, the inequality~\eqref{eq:2.2} means that $\|\Phi\|\le C$. Then, by Theorem~\ref{thm:main_cb}, $\|\Phi\|_{cb}\le C\cdot K$. Put $M=\|\Phi\|_{cb}$. Recall also that without loss of generality we assume that all vectors $\{x_k\}_{k=1}^n$ and $\{y_k\}_{k=1}^n$ are nonzero.
	
	We apply Lemma~\ref{lem:Stinespring} to the map $\Psi=M^{-1}\cdot\Phi$ and obtain a Hilbert space $\K$, isometries $V_{1,2}:H\to\K$ and a $*$-homomorphism $\pi:\ell^\infty_n\to B(\K)$ such that 
	\begin{equation}\label{eq:Stinespring}
		\Phi(a)=M\cdot V_1^* \pi(a)V_2.
	\end{equation}
	Since $\pi$ is a $*$-homomorphism, for every $1\le k\le n$ the operator $\pi(e_k)\in B(\K)$ is self-adjoint and idempotent, hence it is an orthogonal projection on a closed subspace $L_k\subset\K$: $\pi(e_k)=P_{L_k}$. Moreover, again using the fact that $\pi$ is a homomorphism, we conclude that for $k\neq l$ we have $\pi(e_k)\pi(e_l) = 0$, which means that the subspaces $\{L_k\}_{k=1}^n$ are pairwise orthogonal.
	
	Next, since $V_2$ is an isometry, we can identify $H$ and $V_2(H)$, i.e. without loss of generality we may assume that $H\subset \K$ and $V_2=\mathrm{Id}$. Then, since $V_1:H\to\K$ is an isometry, there exists a unitary operator $U:\K\to\K$ such that $V_1=U|_{H}$. The adjoint operator $V_1^*$ then has the form
	\begin{equation}
		V_1^*=P_HU^{-1}.
	\end{equation}
	
	Summing up, if we apply the formula~\eqref{eq:Stinespring} to $a=e_k$, we get for every $u\in H$
	\begin{equation}
		\Phi(e_k)u = M\cdot P_HU^{-1}P_{L_k}u.
	\end{equation}
	
	Put $\widetilde{L}_k = P_{L_k}H\subset L_k$. Then
	\begin{equation}
		\Phi(e_k)u = M\cdot P_HU^{-1}P_{\widetilde{L}_k}u.
	\end{equation}
	On the other hand, recall that according to the formula~\eqref{eq:4.2} we have $\Phi(e_k)u = \la u, y_k\ra x_k$, so for every $u\in H$ we can write
	\begin{equation}\label{eq:4.7}
		\la u, y_k\ra x_k = M\cdot P_HU^{-1}P_{\widetilde{L}_k}u.
	\end{equation}
	
	\subsection{}
	
	The left-hand side of the formula~\eqref{eq:4.7} is a rank one operator. We apply this formula for all $u\in H$ and conclude that $(P_HU^{-1})(\widetilde{L}_k)$ is a one-dimensional subspace of $H$ (which consists of vectors collinear to $x_k$). Then for every $1\le k\le n$ we can choose an orthonormal basis $\{\xi_l^{(k)}\}_l$ in $\widetilde{L}_k$ such that for $l\neq 1$ the vector $U^{-1}\xi_l^{(k)}$ is orthogonal to $x_k$ (this condition is equivalent to the fact that for $l\neq 1$ we have $\la \xi_l^{(k)}, P_{\widetilde{L}_k} Ux_k \ra = 0$, i.e., for $l\neq 1$ the vectors $\xi_l^{(k)}$ should be orthogonal to a fixed vector from $\widetilde{L}_k$, and we can clearly choose an orthonormal basis which satisfies this property).
	
	The orthogonal projection onto $\widetilde{L}_k$ can be written in the coordinates which we have just chosen:
	\begin{equation}
		P_{\widetilde{L}_k} u =\sum_l \la u,\xi_l^{(k)}\ra \xi_l^{(k)},
	\end{equation}
	and therefore we can rewrite the identity \eqref{eq:4.7} in the following way:
	\begin{equation}
		\la u, y_k\ra x_k = M\cdot \sum_{l} \la u, \xi_l^{(k)}\ra P_H U^{-1}\xi_l^{(k)}.
	\end{equation}
	The summands corresponding to $l\neq 1$ in the right-hand side of this formula are orthogonal to $x_k$, hence in fact we have
	\begin{equation}\label{eq:4.10}
		\la u, y_k\ra x_k = M\cdot \la u, \xi_1^{(k)}\ra P_H(U^{-1}\xi_1^{(k)}) =  \la u, M^{1/2}P_H\xi_1^{(k)}\ra M^{1/2}\cdot P_H(U^{-1}\xi_1^{(k)}).
	\end{equation}
	This equation means that the vectors $x_k$ and $P_H(U^{-1}\xi_1^{(k)})$ are collinear (and since $x_k, y_k\neq 0$, both of them are non-zero), hence we can choose nonzero $\alpha_k\in\C$ so that $\alpha_k x_k = M^{1/2}\cdot P_H(U^{-1}\xi_1^{(k)})$. Then by the formula~\eqref{eq:4.10} we have $\overline{\alpha}_k^{-1}y_k = M^{1/2}P_H\xi_1^{(k)}$. Then, since $\{\xi_1^{(k)}\}_{k=1}^n$ is an orthonormal system in $\K$, for every $u\in H$ we have
	\begin{equation}
		\sum_{k=1}^n |\la u, \alpha_k x_k\ra|^2 = M\cdot \sum_{k=1}^n |\la u, U^{-1}\xi_1^{(k)}\ra|^2 \le M \|u\|^2,
	\end{equation}
	which means that the sequence $\{\alpha_kx_k\}_{k=1}^n$ is $M$-Bessel. A similar computation also shows that $\{\overline{\alpha}_k^{-1}y_k\}_{k=1}^n$ is $M$-Bessel. Since $M=\|\Phi\|_{cb}\le C\cdot K$, this choice of the numbers $\{\alpha_k\}_{k=1}^n$ proves Theorem~\ref{thm:main_2}.


\end{document}